# STOCHASTIC DOMINATION FOR A HIDDEN MARKOV CHAIN WITH APPLICATIONS TO THE CONTACT PROCESS IN A RANDOMLY EVOLVING ENVIRONMENT

By Erik I. Broman

*Chalmers University of Technology*

The ordinary contact process is used to model the spread of a disease in a population. In this model, each infected individual waits an exponentially distributed time with parameter 1 before becoming healthy. In this paper, we introduce and study the contact process in a randomly evolving environment. Here we associate to every individual an independent two-state, $\{0,1\}$, background process. Given $\delta_0 < \delta_1$, if the background process is in state 0, the individual (if infected) becomes healthy at rate $\delta_0$, while if the background process is in state 1, it becomes healthy at rate $\delta_1$. By stochastically comparing the contact process in a randomly evolving environment to the ordinary contact process, we will investigate matters of extinction and that of weak and strong survival. A key step in our analysis is to obtain stochastic domination results between certain point processes. We do this by starting out in a discrete setting and then taking continuous time limits.

**1. Introduction.** The first part of this introduction will discuss the concept of stochastic domination and then move on to state our discrete time results. We will then proceed by defining the contact process in a randomly evolving environment (CPREE) that we introduce in this paper. We would like to point out that a model called the contact process in a random environment (CPRE) has been studied before. The first papers concerning this latter model were [2] and [8], and then further studies were carried out in, for instance, [1, 9, 13] and [14]. However, the random environments in those papers are static while here they change over time.

In this paper we are concerned with models on connected graphs $G = (S, E)$ of bounded degree, in which every site $s \in S$ can take values 0 or 1. Here $\sigma$ and $\xi$ will mainly denote configurations on $S$, that is, $\sigma, \xi \in \{0,1\}^S$.

We say that $\xi \preceq \tilde{\xi}$ if $\xi(s) \leq \tilde{\xi}(s)$ for every $s \in S$. An increasing function $f$ is a function $f:\{0,1\}^S \to \mathbb{R}$ such that $f(\xi) \leq f(\tilde{\xi})$ for all $\xi \preceq \tilde{\xi}$. For two probability measures $\mu, \mu'$ on $\{0,1\}^S$, we write $\mu \preceq \mu'$ if for every continuous increasing function $f$ we have that $\mu(f) \leq \mu'(f)$. [$\mu(f)$ is shorthand for $\int f(x)\, d\mu(x)$.] Strassens theorem (see [10], page 72) states that if $\mu \preceq \mu'$, then there exist random variables $X, X'$ with distributions $\mu, \mu'$, respectively, defined on the same probability space, such that $X \preceq X'$ a.s.

We will need the following standard definition.

Definition 1.1. Let $S$ be such that $|S| < \infty$ and let $\mu$ be a probability measure on $\{0,1\}^S$ with full support. $\mu$ is said to be monotone, if for every $s \in S$ and any $\xi, \tilde{\xi} \in \{0,1\}^{S \setminus s}$ such that $\xi \preceq \tilde{\xi}$, one has that

$$\mu(\sigma(s) = 1 | \sigma(S \setminus s) \equiv \xi) \leq \mu(\sigma(s) = 1 | \sigma(S \setminus s) \equiv \tilde{\xi}).$$

If $|S| = \infty$, we say that a probability measure $\mu$ on $\{0,1\}^S$ is monotone if the restriction of $\mu$ to any finite subset of $S$ is monotone.

For $p \in [0,1]$, let each site $s \in S$, independently of all others, take value 1 with probability $p$ and 0 with probability $1-p$. Write $\pi_p$ for this product measure on $\{0,1\}^S$. For any probability measure $\mu$ on $\{0,1\}^S$ define $p_{\max,\mu}$ by

$$p_{\max,\mu} := \sup\{p \in [0,1] : \pi_p \preceq \mu\}.$$

The supremum is easily seen to be obtained, which motivates the notation. Similarly define

$$p_{\min,\mu} := \inf\{p \in [0,1] : \mu \preceq \pi_p\}.$$

Next, informally here we will think of $\{B_n\}_{n=1}^{\infty}$ as a background process which influences another process $\{X_n\}_{n=1}^{\infty}$. Formally, fix $0 \leq \alpha_0 \leq \alpha_1 \leq 1$ and let $\{B_n\}_{n=1}^{\infty}$ be any process with state space $\{0,1\}$. Conditioned on $\{B_n\}_{n=1}^{\infty}$ let the process $\{X_n\}_{n=1}^{\infty}$, also with state space $\{0,1\}$, be a sequence of conditionally independent random variables where the (conditional) distribution of $X_k$ is

$$\begin{array}{lll} \text{if} & \text{then} & \text{w.p.} \\ B_k = 0 & X_k = 1 & \alpha_0 \\ B_k = 1 & X_k = 1 & \alpha_1 \end{array} \tag{1}$$

for every $k \geq 1$.

We will say that $\mu$ is translation invariant on $\mathbb{N}$ if for every $l \geq 1$, $k \geq 0$ and any $\xi \in \{0,1\}^{\{1,\ldots,l\}}$

$$\mu(\sigma(1,\ldots,l) \equiv \xi) = \mu(\sigma(k+1,\ldots,k+l) \equiv \xi).$$

In Section 2 we will prove the following proposition.

PROPOSITION 1.2. *Assume that the distribution of* $\{B_n\}_{n=1}^{\infty}$ *is monotone and translation invariant. Then the sequence*

$$\{\mathbb{P}(X_n = 1 | X_{n-1} = \cdots = X_1 = 0)\}_{n \geq 1},$$

*is decreasing in* $n$*. In addition the limit equals* $p_{\max,\mu}$*, where* $\mu$ *is the distribution of* $\{X_n\}_{n=1}^{\infty}$*.*

The proof is an easy consequence of results from [5] and [12].

We are now ready to define the discrete background process that we will use throughout this paper. For $p, \gamma \in [0,1]$, define the Markov chain $\{B_n\}_{n=1}^{\infty}$ in the following way:

$$B_1 = \begin{cases} 1, & \text{w.p. } p, \\ 0, & \text{w.p. } 1-p, \end{cases} \tag{2}$$

and for $k \geq 2$,

$$\begin{array}{lll} \text{if} & \text{then} & \text{w.p.} \\ B_{k-1} = 0 & B_k = 1 & \gamma p \\ B_{k-1} = 1 & B_k = 0 & \gamma(1-p). \end{array} \tag{3}$$

In other words, $B_k$ takes the same value as $B_{k-1}$ unless there is an update which happens independently with probability $\gamma$. If an update occurs, $B_k = 1$ with probability $p$, and $B_k = 0$ with probability $1-p$. Using (1) this defines a joint process $\{(B_n, X_n)\}_{n=1}^{\infty}$, whose second marginal is an example of a so-called hidden Markov chain. The main theorem of Section 3 is the following, here $\mu$ refers to the distribution of $\{X_n\}_{n=1}^{\infty}$ with the background process as above.

THEOREM 1.3. *We have that*

$$p_{\max,\mu} = \tfrac{1}{2}(1 - C - \sqrt{(1-C)^2 - 4D}),$$

*where*

$$C = (1 - \alpha_0 - \alpha_1) - \gamma(1 - \alpha_0 - (1-p)(\alpha_1 - \alpha_0))$$

*and*

$$D = \alpha_0\alpha_1 + \gamma(\alpha_1(1-\alpha_0) - (1-p)(\alpha_1 - \alpha_0)).$$

*Furthermore,*

$$p_{\min,\mu} = \tfrac{1}{2}(1 + C' + \sqrt{(1-C')^2 - 4D'}),$$

*where* $C'$ *and* $D'$ *are as* $C$ *and* $D$ *but with* $\alpha_0, \alpha_1, p$ *and* $\gamma$ *replaced by* $1-\alpha_1$*,* $1-\alpha_0$*,* $1-p$ *and* $\gamma$*, respectively.*

REMARK. It is easy to check that $p_{\max,\mu}$ ($p_{\min,\mu}$) is increasing (decreasing) in $\gamma$. This is natural since intuitively, as $\gamma$ increases, $\{X_n\}_{n=1}^{\infty}$ looks more like an i.i.d. process and so our approximations should get better.

The proof of Theorem 1.3 unfortunately involves some tedious (but straightforward) calculations; however this result is needed for all the other results of this paper.

From the results of Section 3, we will in Section 4 prove our next result. First define

$$X_n^c := \sum_{i=1}^{n} X_i \qquad \forall n \in \mathbb{N},$$

where $c$ indicates that we are counting the number of 1's up to time $n$. The pair of processes $\{(B_t, X_t)\}_{t\geq 0}$, to be defined below will be a Markov process. Furthermore, it will be the continuous time analogue of the pair of processes $\{(B_n, X_n^c)\}_{n=1}^{\infty}$. To define $\{(B_t, X_t)\}_{t\geq 0}$, let $B_0 = 1$ with probability $p$ and $B_0 = 0$ with probability $1-p$, also let $X_0 = 0$. We define the transition rates at time $t \geq 0$, for the Markov process $\{(B_t, X_t)\}_{t\geq 0}$ as follows:

| from | to | with intensity |
|---|---|---|
| $(0,k)$ | $(1,k)$ | $\gamma p$ |
| $(1,k)$ | $(0,k)$ | $\gamma(1-p)$ |
| $(0,k)$ | $(0,k+1)$ | $\alpha_0$ |
| $(1,k)$ | $(1,k+1)$ | $\alpha_1$ |

for any $k \geq 0$. Informally this can be described in the following way. Letting the value of $B_0$ be chosen as before, the $\{B_t\}_{t\geq 0}$ process waits an exponentially distributed time with parameter $\gamma > 0$ before it updates its status. After an update, this process takes value 1 with probability $p$ and 0 with probability $1-p$, and all of this is done independently of everything else. If the background process is in state 0, $\{X_t\}_{t\geq 0}$ increases by one every time a Poisson process with rate $\alpha_0$ has an arrival. If instead the background process is in state 1, $\{X_t\}_{t\geq 0}$ increases by one every time a Poisson process with rate $\alpha_1$ has an arrival. In short, $\{X_t\}_{t\geq 0}$ is the counting process for a type of Poisson process where the parameter comes from $\{B_t\}_{t\geq 0}$. Analogous to the definition of $p_{\max,\mu}$, define $\lambda_{\max,\mu}$, where $\mu$ here refers to the distribution of $\{X_t\}_{t\geq 0}$, in the following way. $\lambda_{\max,\mu}$ is the maximum real number $\lambda$ such that a Poisson($\lambda$)-process can be coupled with the process $\{X_t\}_{t\geq 0}$ so that if the Poisson($\lambda$)-process has an arrival at time $\tau \in [0, \infty)$ then so does the $\{X_t\}_{t\geq 0}$ process. In other words, there exists $\{X_t'\}_{t\geq 0}$ with distribution Poisson($\lambda$) coupled with $\{X_t\}_{t\geq 0}$ such that

$$X_t - X_t' \text{ is nondecreasing in } t.$$

Define $\lambda_{\min,\mu}$ to be the minimal real number $\lambda$ such that a Poisson($\lambda$)-process can be coupled with the process $\{X_t\}_{t\geq 0}$ so that if the $\{X_t\}_{t\geq 0}$ process has an arrival at time $\tau \in [0,\infty)$ then so does the Poisson($\lambda$)-process. Observe that $\lambda_{\max,\mu} = \lambda_{\max,\mu}(\alpha_0,\alpha_1,\gamma,p)$. We will write out the arguments in most equations, but not in more general discussions. Trivially $\alpha_0 \leq \lambda_{\max,\mu} \leq \lambda_{\min,\mu} \leq \alpha_1$. For future convenience let $\mathrm{Poi}^{\gamma,p}_{\alpha_0,\alpha_1}$ denote the distribution of $\{X_t\}_{t\geq 0}$ and $\mathrm{Poi}_\lambda$ denote the distribution of a Poisson process with intensity $\lambda$. The coupling described above is a form of stochastic domination and we will write

$$\mathrm{Poi}_{\lambda_{\max,\mu}} \preceq \mathrm{Poi}^{\gamma,p}_{\alpha_0,\alpha_1}, \tag{4}$$

and

$$\mathrm{Poi}^{\gamma,p}_{\alpha_0,\alpha_1} \preceq \mathrm{Poi}_{\lambda_{\min,\mu}}. \tag{5}$$

Define

$$\begin{aligned}\bar\lambda &= \bar\lambda(\alpha_0,\alpha_1,\gamma,p)\\ &:= \tfrac{1}{2}\left(\alpha_0+\alpha_1+\gamma-\sqrt{(\alpha_1-\alpha_0-\gamma)^2+4\gamma(1-p)(\alpha_1-\alpha_0)}\right).\end{aligned}$$

Before stating our next result we would like to point out that the pair of processes $\{(B_t,X_t)\}_{t\geq 0}$ has been studied before (see [7] and some of the references therein). However, in these papers the focus and motivation for the study is completely different from ours, the question of interest being how to best determine $B_0$ by observing the $\{X_t\}_{t\geq 0}$-process. See further the remark after Theorem 1.5 below.

THEOREM 1.4. *Let $\{(B_t,X_t)\}_{t\geq 0}$ be as above. For every choice of $\alpha_0$, $\alpha_1$, $\gamma > 0$ with $\alpha_0 \leq \alpha_1$ and $p \in [0,1]$ we have that*

$$\lambda_{\max,\mu}(\alpha_0,\alpha_1,\gamma,p) = \bar\lambda, \tag{6}$$

*and for $p > 0$*

$$\lambda_{\min,\mu}(\alpha_0,\alpha_1,\gamma,p) = \alpha_1.$$

REMARKS.

- Note the apparent lack of symmetry between $\lambda_{\max,\mu}$ and $\lambda_{\min,\mu}$. Informally, consider for a moment the model to be a point process, where the process is 0 unless there is an arrival, in which case it takes the value 1. We can then see that the true symmetric statement of the $\lambda_{\max,\mu}$ result would concern a model which is 1 unless there is an arrival in which case it takes the value 0. This however does not correspond to the result concerning $\lambda_{\min,\mu}$.

- We will show in Section 4 that $\lambda_{\max,\mu}(\alpha_0,\alpha_1,\gamma,p)\to\min(\alpha_1,\alpha_0+\gamma)$ as $p\to 1$. Hence if $\gamma>\alpha_1-\alpha_0$, then $\lambda_{\max,\mu}(\alpha_0,\alpha_1,\gamma,p)\to\alpha_1$ as $p\to 1$ which one would expect. In contrast, for every $p>0$ $\lambda_{\min,\mu}(\alpha_0,\alpha_1,\gamma,p)=\alpha_1$ and so $\lambda_{\min,\mu}(\alpha_0,\alpha_1,\gamma,p)\not\to\alpha_0$ as $p\to 0$ as one might have suspected; this gives a discontinuity at $p=0$. Also, it is trivial to show that $\lambda_{\max,\mu}(\alpha_0,\alpha_1,\gamma,p)\to\alpha_0$ as $p\to 0$.
- Intuitively, as $\gamma$ grows larger, the suppressing of possible arrivals becomes "increasingly independent." Whenever a possible arrival occurs, the background process has with very high probability been updated since the last possible arrival, and if so the new arrival is suppressed independently of everything else. Therefore, as $\gamma$ grows larger, we would expect our process to look more and more like an ordinary Poisson process. As a consequence, we would expect our approximation to get better. This is confirmed by studying the derivative of $\bar{\lambda}$ with respect to $\gamma$; it is easy to check that $\bar{\lambda}$ is increasing in $\gamma$. Furthermore, for fixed $0<p<1$, it follows from the proof of Proposition 1.9, where we take the limit $\gamma\to\infty$ in equation (6) above, that $\lim_{\gamma\to\infty}\lambda_{\max,\mu}(\alpha_0,\alpha_1,\gamma,p)=p\alpha_1+(1-p)\alpha_0$. Of course this is exactly what you would expect to get. Also, by letting $\gamma\to 0$ in equation (6) we get that $\lambda_{\max,\mu}(\alpha_0,\alpha_1,0,p)\to\alpha_0$. This last result is also natural. As $\gamma$ becomes smaller, we will find longer and longer time intervals in which the background process is in the lower state. Therefore the Poisson process we dominate must have lower and lower density.

It is natural to ask for quantitative versions of Theorem 1.4 for finite time, and in fact we will use such results to prove Theorem 1.4. Therefore, for $T>0$, let $\lambda^T_{\max,\mu}(\alpha_0,\alpha_1,\gamma,p)$ denote the maximum intensity of the Poisson process which the second marginal of the truncated process $\{(B_t,X_t)\}_{t\in[0,T]}$ can dominate. Define $\lambda^T_{\min,\mu}(\alpha_0,\alpha_1,\gamma,p)$ analogously. We feel that this bound is interesting in its own right and we therefore present it in our next theorem together with a lower bound on $\lambda^T_{\max,\mu}(\alpha_0,\alpha_1,\gamma,p)$ and a result for $\lambda^T_{\min,\mu}(\alpha_0,\alpha_1,\gamma,p)$ (this last result will follow from the proof of Theorem 1.4).

Let

$$L=L(\alpha_0,\alpha_1,\gamma,p):=\sqrt{(\alpha_1-\alpha_0-\gamma)^2+4\gamma(1-p)(\alpha_1-\alpha_0)}.$$

THEOREM 1.5. *For every choice of $\alpha_0,\alpha_1,\gamma,T>0$ with $\alpha_0\le\alpha_1$ and $p\in(0,1)$ we have that*

$$\lambda^T_{\max,\mu}(\alpha_0,\alpha_1,\gamma,p)\ge\bar{\lambda}+(p\alpha_1+(1-p)\alpha_0-\bar{\lambda})e^{-TL}. \tag{7}$$

*Furthermore there exists a constant $E > 0$, depending on $\alpha_1, \alpha_0, \gamma$ and $p$, such that*

$$\lambda^T_{\max,\mu}(\alpha_0,\alpha_1,\gamma,p) \le \bar\lambda + \frac{1}{T}(p\alpha_1 + (1-p)\alpha_0 - \bar\lambda)\frac{1-e^{-TE}}{E}. \tag{8}$$

*Finally*

$$\lambda^T_{\min,\mu} = \lambda_{\min,\mu} = \alpha_1.$$

REMARKS.

- Observe that the right-hand side of equation (7) tends to $p\alpha_1 + (1-p)\alpha_0$ as $T$ tends to 0, and that it tends to $\bar\lambda = \lambda_{\max,\mu}(\alpha_0,\alpha_1,\gamma,p)$ as $T$ tends to infinity. Both results are of course what you would expect. The same is true for the upper bound of equation (8).
- In [7] they calculate the probability density $F_{ij}(t)$ of the first arrival of $\{X_t\}_{t\ge 0}$ occurring at time $t$ while $B_t = j$ given that $B_0 = i$. Through a nontrivial although straightforward series of calculations it is possible from their results to arrive at the conclusion that

  $$\lambda^T_{\max,\mu}(\alpha_0,\alpha_1,\gamma,p) \le \bar\lambda + \frac{1}{T}c,$$

  for some constant $c$. However, this result does not give the right asymptotic behavior as $T$ tends to 0. Furthermore, since this approach also rests on work needed to arrive at the expression for the probability density $F_{ij}(t)$ we do not use the results of [7] in our proof.

### 1.1. *Models.*

1.1.1. *The contact process.* In this section we will discuss some basic concepts concerning the contact process, see [10] for results up to 1985 and [11] for results between 1985 and 1999. Consider a graph $G = (S, E)$ of bounded degree. In the contact process the state space is $\{0,1\}^S$. We will let 1 represent an infected individual, while a 0 will be used to represent a healthy individual. Let $\lambda > 0$, and define the flip rate intensities to be

$$C(s,\sigma) = \begin{cases} 1, & \text{if } \sigma(s) = 1, \\ \lambda \sum_{(s',s)\in E} \sigma(s'), & \text{if } \sigma(s) = 0. \end{cases} \tag{9}$$

By flip rate intensities, informally, we mean as usual that every site $s \in S$ waits an exponentially distributed time with parameter $C(s,\sigma)$ before changing its status. Here, $\mathbf{1}_0, \mathbf{1}_1$ will denote the measures that put mass one on the configuration of all 0's and all 1's, respectively. If we let the initial distribution be $\sigma \equiv 1$, the distribution of this process at time $t$, which we will denote by $\mathbf{1}_1 T_\lambda(t)$, is known to converge as $t$ tends to infinity. This is simply because it is a so-called "attractive" process and $\sigma \equiv 1$ is the maximal state; see [10], page 265. This limiting distribution will be referred to as the upper invariant measure for the contact process with parameter $\lambda$ and will be denoted by $\nu_\lambda$. We then let $\Psi^\lambda$ denote the stationary Markov process on $\{0,1\}^S$ with initial (and invariant) distribution $\nu_\lambda$. One can also choose to start the process with any set $A \subset S$ of infected individuals and then use the flip rate intensities above. Denote this latter process by $\Psi^{\lambda,A}$. We say that the process dies out if for any $s \in S$

$$\Psi^{\lambda,\{s\}}(\sigma_t \not\equiv 0\ \forall t \geq 0) = 0,$$

and otherwise it survives. We also say that the process survives strongly if

$$\Psi^{\lambda,\{s\}}(\sigma_t(s) = 1 \text{ i.o.}) > 0.$$

We say that the process survives weakly if it survives but does not survive strongly. These and all other statements like it, made here and later, are independent of the specific choice of the site $s$; see [11]. We will use the same definition of survival for some closely related processes below. It is well known that for any graph (see [11], page 42) there exists two critical parameter values $0 \leq \lambda_{c1} \leq \lambda_{c2} \leq \infty$ such that:

- $\Psi^{\lambda,\{s\}}$ dies out if $\lambda < \lambda_{c1}$,
- $\Psi^{\lambda,\{s\}}$ survives weakly if $\lambda_{c1} < \lambda < \lambda_{c2}$,
- $\Psi^{\lambda,\{s\}}$ survives strongly if $\lambda > \lambda_{c2}$.

The above description of the contact process with flip rate intensities chosen according to (9) is standard. However for our purposes it is more convenient to use the following flip rate intensities. Let $\delta > 0$ and

$$C(s,\sigma) = \begin{cases} \delta, & \text{if } \sigma(s) = 1, \\ \displaystyle\sum_{(s',s)\in E} \sigma(s'), & \text{if } \sigma(s) = 0. \end{cases} \tag{10}$$

This is just a time scaling of the original model. We will denote the upper invariant measure by $\nu_\delta$ and the corresponding process starting with distribution $\nu_\delta$ by $\Psi_\delta$. If we instead choose to start with a specific set $A \subset S$ of infected individuals, we denote the corresponding process by $\Psi_\delta^A$. We will let the distribution of the process at time $t \geq 0$ be denoted by $\nu_{\delta,t}^A$. At some point we need to consider the process $\Psi_\delta^{\lambda,A}$, this is exactly like the model just described except for a $\lambda$ inserted in front of the sum in equation (10).

As above, it follows that there exists $0 \leq \delta_{c1} \leq \delta_{c2} \leq \infty$ such that:

- $\Psi_\delta^{\{s\}}$ dies out if $\delta > \delta_{c2}$,
- $\Psi_\delta^{\{s\}}$ survives weakly if $\delta_{c1} < \delta < \delta_{c2}$,

- $\Psi_\delta^{\{s\}}$ survives strongly if $\delta < \delta_{c1}$.

We point out that on $G = \mathbb{Z}^d$ it is known (see [3]) that $\delta_{c1} = \delta_{c2}$. It is also well known (see [11]) that on any homogeneous tree of degree larger than or equal to 3, this is not the case.

1.1.2. *CPREE.* This model is a pair of processes $\{(B_t, Y_t)\}_{t\geq 0}$ with state space $\{\{0,1\} \times \{0,1\}\}^S$. The second coordinate of $\{(B_t, Y_t)\}_{t\geq 0}$ will represent whether an individual is infected or not, while the first coordinate will represent how prone the individual is to recover. With a slight abuse of notation we have chosen the first coordinate to be denoted by $\{B_t\}_{t\geq 0}$ even though a process with this notation was already defined previously in the Introduction. However, at every site $s \in S$, the marginal of the $\{B_t\}_{t\geq 0}$ process defined in this subsection [denoted by $\{B_t(s)\}_{t\geq 0}$] will be independent of the rest of the $\{B_t\}_{t\geq 0}$ process defined here, and have distribution according to the process with the same notation defined earlier. It will be clear from context which of these two we are referring to.

For any $A \subset S$, let $Y_0(s) = 1$ iff $s \in A$, and let $B_0 \sim \pi_p$. For $0 \leq \delta_0 < \delta_1 < \infty$, $\gamma > 0$ and $p \in [0,1]$, let the flip rate intensities $C(s, (B_t, Y_t))$ of a site $s \in S$ be

| from | to | with intensity |
|---|---|---|
| $(0,0)$ | $(1,0)$ | $\gamma p$ |
| $(0,1)$ | $(1,1)$ | $\gamma p$ |
| $(1,0)$ | $(0,0)$ | $\gamma(1-p)$ |
| $(1,1)$ | $(0,1)$ | $\gamma(1-p)$ |
| $(0,0)$ | $(0,1)$ | $\sum_{(s',s)\in E} Y_t(s')$ |
| $(1,0)$ | $(1,1)$ | $\sum_{(s',s)\in E} Y_t(s')$ |
| $(0,1)$ | $(0,0)$ | $\delta_0$ |
| $(1,1)$ | $(1,0)$ | $\delta_1$. |

Denote the distribution of $\{Y_t\}_{t\geq 0}$ by $\Psi_{\delta_0,\delta_1}^{\gamma,p,A}$ and the distribution at a fixed time $t$ by $\nu_{\delta_0,\delta_1,t}^{\gamma,p,A}$. The definition of dying out, surviving weakly and surviving strongly is the same as for the ordinary contact process. At this point a question naturally arises. For fixed $\delta_0, \delta_1, \gamma$ and $p$, do the initial state of the background process have any effect on this definition? This point is raised in Section 6, where we list some open questions. Note that we are here assuming that $B_0 \sim \pi_p$ which then is included in the definition.

We will write

$$\Psi_{\delta_0,\delta_1}^{\gamma,p,A} \preceq \Psi_\delta^A$$

when we mean that there exists a process $\{Y_t\}_{t\geq 0}$ with distribution as above and a process $\{Y'_t\}_{t\geq 0}$ with distribution $\Psi_\delta^A$ coupled such that

$$Y_t(s) \leq Y'_t(s) \qquad \forall s \in S \text{ and } \forall t \geq 0,$$

and use the obvious notation for all similar types of situations. This stochastic ordering also implies that

$$\nu_{\delta_0,\delta_1,t}^{\gamma,p,A} \preceq \nu_{\delta,t}^A \qquad \forall t \geq 0.$$

It is easy to show that $\Psi_{\delta_0,\delta_1}^{\gamma,p,A}$ is in this sense stochastically decreasing in $p$. We have already introduced this notation for continuous time processes in (4) and (5). There it was a relation between jump processes indexed by $t \geq 0$, while here it is a relation between processes with state space $\{0,1\}^S$. It will be clear from the context which one we are referring to. By the recovery process of a site, we mean the process governing the recoveries of that site. In the ordinary contact process it is a $\mathrm{Poi}_\delta$ process, while by the definition of the CPREE above, the recovery process at every site is in fact a $\mathrm{Poi}_{\delta_0,\delta 1}^{\gamma,p}$ process as defined earlier. This explains the relation between Theorem 1.4 and our next result.

We will let $\Delta_G$ denote the maximum degree of a graph $G$ of bounded degree. We can now list our main results concerning this model:

THEOREM 1.6. *Let $G = (S,E)$ be any graph of bounded degree and $A \subset S$, be such that $|A| < \infty$. For any $\delta < \min(\delta_1, \delta_0 + \gamma)$ there exists a $p = p(\delta, \delta_0, \delta_1, \gamma) \in (0,1)$ large enough so that*

$$\Psi_{\delta_0,\delta_1}^{\gamma,p,A} \preceq \Psi_\delta^A. \tag{11}$$

*Furthermore, for $\delta > \min(\delta_1, \delta_0 + \gamma)$ there is no $p \in (0,1)$ such that (11) holds.*

Our next result also uses Theorem 1.4. However, it does so in a different way. The reason for this is that a straightforward approach would need a result for $\lambda_{\min,\mu}(\alpha_0,\alpha_1,\gamma,p)$ analogous to the one we have for $\lambda_{\max,\mu}(\alpha_0,\alpha_1,\gamma,p)$. However this is false since $\lambda_{\min,\mu}(\alpha_0,\alpha_1,\gamma,p)$ is equal to $\alpha_0$ for any $p > 0$. Here, let $\Psi_{\delta_0,\delta_1,B_0(A)\equiv 0}^{\gamma,p,A}$ denote the distribution of $\Psi_{\delta_0,\delta_1}^{\gamma,p,A}$ conditioned on the event that $B_0(s) = 0$ for every $s \in A$.

THEOREM 1.7. *Let $G = (S,E)$ be any graph of bounded degree. Let $A \subset S$, be such that $|A| < \infty$ and $\gamma \geq \Delta_G$. For any choice of $\delta > \delta_0$ and $\lambda < 1$ there exists a $p = p(\delta, \lambda, \delta_0, \delta_1, \gamma) \in (0,1)$ small enough so that*

$$\Psi_\delta^{\lambda,A} \preceq \Psi_{\delta_0,\delta_1,B_0(A)\equiv 0}^{\gamma,p,A}.$$

REMARKS. It is unfortunate that we need the assumption that $B_0(s)=0$ for every $s\in A$. However, this is of no importance when we later apply the theorem to prove Theorem 1.8 stated below. Also, as $\lambda\to 1$, the proof requires that $p\to 0$. The hypothesis $\gamma\ge\Delta_G$ and $\lambda<1$, might look artificial, however a lower bound on $\gamma$ is required as Example 5.1 shows, and for $\lambda=1$, the statement is false as Example 5.2 shows. We thank the referee for pointing out these two examples.

We are now ready to state the main theorem concerning the CPREE model of this paper. Results 1–3 use an easy coupling argument while 4–6 follow from applications of Theorems 1.6 and 1.7. Here, any statements similar to $p_{c1}<p<p_{c2}$ in the case that $p_{c1}=p_{c2}$ should be interpreted as empty statements.

THEOREM 1.8. *Let $s\in S$, $0\le\delta_0\le\delta_1<\infty$ and consider the process $\Psi^{\gamma,p,\{s\}}_{\delta_0,\delta_1}$. We have the following results:*

1. *Assume that $\delta_{c1}<\delta_0<\delta_{c2}<\delta_1$. There exists $p_{c2}=p_{c2}(\delta_0,\delta_1,\gamma)\in[0,1]$ such that $\Psi^{\gamma,p,\{s\}}_{\delta_0,\delta_1}$ dies out if $p>p_{c2}$ and survives weakly if $p<p_{c2}$.*
2. *Assume that $\delta_0<\delta_{c1}\le\delta_{c2}<\delta_1$. There exists $p_{c2}=p_{c2}(\delta_0,\delta_1,\gamma)\in[0,1]$ and $p_{c1}=p_{c1}(\delta_0,\delta_1,\gamma)\in[0,1]$ such that $p_{c1}\le p_{c2}$ and $\Psi^{\gamma,p,\{s\}}_{\delta_0,\delta_1}$ dies out if $p>p_{c2}$ survives weakly if $p_{c1}<p<p_{c2}$ and survives strongly if $p<p_{c1}$.*
3. *Assume that $\delta_0<\delta_{c1}<\delta_1<\delta_{c2}$. There exists $p_{c1}=p_{c1}(\delta_0,\delta_1,\gamma)\in[0,1]$ such that $\Psi^{\gamma,p,\{s\}}_{\delta_0,\delta_1}$ survives strongly if $p<p_{c1}$ and survives weakly if $p>p_{c1}$.*
4. *In case number 1, if $\gamma>\delta_{c2}-\delta_0$ then $p_{c2}<1$ and if $\gamma\ge\Delta_G$, then $p_{c2}>0$.*
5. *In case number 2, if $\gamma>\delta_{c2}-\delta_0$ then $p_{c2}<1$, if $\gamma>\delta_{c1}-\delta_0$ then $p_{c1}<1$ and if $\gamma\ge\Delta_G$, then $p_{c1},p_{c2}>0$.*
6. *In case number 3, if $\gamma>\delta_{c1}-\delta_0$ then $p_{c1}<1$ and if $\gamma\ge\Delta_G$, then $p_{c1}>0$.*

REMARKS.

- We do not include trivial cases like $\delta_{c1}<\delta_0<\delta_1<\delta_{c2}$ in the statement.
- One might suspect that the condition $\gamma>\delta_{c2}-\delta_0$ should in fact be $\gamma>\delta_1-\delta_0$, considering the statement of Theorem 1.6. The point is however that we must only be able to choose the $\delta$ of Theorem 1.6 to be larger than $\delta_{c2}$, not $\delta_1$.
- We would like to point out that even if we only show that $p_{c1},p_{c2}<1$ whenever $\gamma>\delta_{c1}-\delta_0$, $\gamma>\delta_{c2}-\delta_0$, respectively, there is no apparent reason why this should not be true for all $\gamma>0$. Similarly for $p_{c1},p_{c2}>0$.

The rest of the paper is organized as follows. Proposition 1.2 is proved in Section 2. This is then used to prove Theorem 1.3 in Section 3. In Section 4 we use a limiting argument to conclude Theorem 1.4 from Theorem 1.3. We then exploit Theorem 1.4 to prove Theorems 1.6 and 1.7 in Section 5. Finally, these last two theorems will be used to prove Theorem 1.8 in the same section.

We exploit the techniques for proving Theorem 1.6 further to conclude the following results concerning $p_{c1}$ and $p_{c2}$.

PROPOSITION 1.9. *Fix $i \in \{1,2\}$ and assume that $\delta_0 < \delta_{ci}$. We have for $p_{ci}$,*

$$\limsup_{\gamma\to\infty} p_{ci}(\delta_0,\delta_1,\gamma) \le \frac{\delta_{ci}-\delta_0}{\delta_1-\delta_0}.$$

REMARK. We conjecture that the limit exists and that

$$\lim_{\gamma\to\infty} p_{ci}(\delta_0,\delta_1,\gamma) = \frac{\delta_{ci}-\delta_0}{\delta_1-\delta_0}.$$

We cannot prove this with the techniques of this paper; this is closely related to the remarks after Theorem 1.4. However, the intuition why it should be true is that as $\gamma\to\infty$, the recovery process should become increasingly similar to an ordinary Poisson($\delta_0 + p(\delta_1-\delta_0)$)-process (see the remarks of Theorem 1.4). In turn, our CPREE then should become more and more like an ordinary contact process with recovery rate $\delta_0 + p(\delta_1-\delta_0)$. Therefore we should get that $p_c$ solves the equation $\delta_c = \delta_0 + p(\delta_1-\delta_0)$.

As $\gamma$ tends to 0 we can unfortunately not conclude anything about $p_{ci}$. The reason is yet again connected to the remark after Theorem 1.4. We know that for $\gamma = 0$, $\lambda_{\min,\mu}(\alpha_0,\alpha_1,0,p) = \alpha_1$ from Theorem 1.4 and that $\lambda_{\max,\mu}(\alpha_0,\alpha_1,0,p) = \alpha_0$. Therefore the stochastic domination techniques we use in this paper do not yield any nontrivial results. We also point out that the case $\gamma = 0$, corresponds to the CPRE and we therefore refer to the papers mentioned in the first paragraph of the Introduction.

We also have the following easy result about $p_{c1}$, $p_{c2}$.

PROPOSITION 1.10. *We have that for any $\gamma > 0$ and $\delta_1 > \delta_{ci} > \delta_0$, where $i \in \{1,2\}$*

$$\lim_{\delta_0\uparrow\delta_{ci}} p_{ci}(\delta_0,\delta_1,\gamma) = 0.$$

REMARK. One would of course expect that

$$\lim_{\delta_1\downarrow\delta_{ci}} p_{ci}(\delta_0,\delta_1,\gamma) = 1.$$

However, it is not possible to prove this the same way as we prove Proposition 1.10; again this is a fact that propagates from Theorem 1.4.

**2. Proof of Proposition 1.2.** The proof of Proposition 1.2 will require the following two results, the first uses Lemma 3.2 of [5] and the second is a restatement (which is more suitable for our purposes) of Theorem 1.2 of [12].

LEMMA 2.1. *If* $\{B_n\}_{n=1}^{\infty}$ *is monotone then* $\{X_n\}_{n=1}^{\infty}$ *is monotone.*

PROOF. Let $\{Z_n\}_{n=1}^{\infty} \sim \pi_{1-(1-\alpha_1)/(1-\alpha_0)}$ and $\{Z'_n\}_{n=1}^{\infty} \sim \pi_{\alpha_0}$ be independent. Observe that $\{X_n\}_{n=1}^{\infty}$ has the same distribution as $\{\max(\min(B_n, Z_n), Z'_n)\}_{n=1}^{\infty}$. It follows from Lemma 3.2 of [5] that $\{\min(B_n, Z_n)\}_{n=1}^{\infty}$ is monotone. It then follows similarly that $\{\max(\min(B_n, Z_n), Z'_n)\}_{n=1}^{\infty}$ is monotone. □

LEMMA 2.2. *Let* $\mu$ *be a translation invariant measure on* $\{0,1\}^{\mathbb{N}}$ *which is monotone. Then the following two statements are equivalent.*

1. $\pi_p \preceq \mu$.
2. *For any* $n \in \mathbb{N}$

$$\mu(\sigma(n) = 1 | \sigma(1, \dots, n-1) \equiv 0) \geq p.$$

PROOF OF PROPOSITION 1.2. Let

$$A_n := \mathbb{P}(X_n = 1 | X_{n-1} = \cdots = X_1 = 0).$$

Since $\{X_n\}_{n=1}^{\infty}$ is monotone (Lemma 2.1) and translation invariant it is easy to see that $A_n$ is decreasing in n, and therefore the limit $A = \lim_{n\to\infty} A_n$ exists. It is now an easy consequence of Lemma 2.2 that this limit is equal to $p_{\max,\mu}$. □

The above results show that when the assumptions of the theorem hold then

$$\inf_{n\in\mathbb{N}, \xi\in\{0,1\}^{n-1}} \mathbb{P}(X_n = 1 | (X_{n-1}, \dots, X_1) \equiv \xi) = p_{\max,\mu}.$$

It is very easy to find examples for which this statement is not true. For instance let $(X, Y) \in \{0,1\} \times \{0,1\}$ and $\mathbb{P}(X = Y = 1) = \mathbb{P}(X = Y = 0) = 1/2$. This dominates a product measure with positive density so that $p_{\max,\mu} > 0$, while $\mathbb{P}(X = 1 | Y = 0) = 0$.

**3. Discrete time domination results.** This section is devoted to the proof of Theorem 1.3. $\{(B_n, X_n)\}_{n=1}^{\infty}$ are the processes defined in the Introduction. We start with the following lemma; we do not include the elementary proof.

Lemma 3.1. *The process $\{B_n\}_{n=1}^{\infty}$ is monotone.*

We will also need the following lemma which gives us a recursion formula of $A_n$ expressed in terms of $A_{n-1}$.

Lemma 3.2. *We have that*

$$A_n = \frac{CA_{n-1} + D}{1 - A_{n-1}}, \tag{12}$$

*with $C, D$ as in Theorem 1.3.*

Proof. The proof is straightforward; however, it involves some tedious calculations. We have

$$\begin{aligned} A_n &= \frac{\mathbb{P}(X_n = 1, X_{n-1} = 0 | (X_{n-2}, \ldots, X_1) \equiv 0)}{\mathbb{P}(X_{n-1} = 0 | (X_{n-2}, \ldots, X_1) \equiv 0)} \\ &= \frac{\mathbb{P}(X_n = 1, X_{n-1} = 0 | (X_{n-2}, \ldots, X_1) \equiv 0)}{1 - A_{n-1}}. \end{aligned}$$

Observe that

$$\begin{aligned} &\mathbb{P}(X_n = 1 | (X_{n-1}, \ldots, X_1) \equiv 0) \\ &\quad = \mathbb{P}(X_n = 1 | B_n = 1, (X_{n-1}, \ldots, X_1) \equiv 0) \mathbb{P}(B_n = 1 | (X_{n-1}, \ldots, X_1) \equiv 0) \\ &\qquad + \mathbb{P}(X_n = 1 | B_n = 0, (X_{n-1}, \ldots, X_1) \equiv 0) \mathbb{P}(B_n = 0 | (X_{n-1}, \ldots, X_1) \equiv 0) \\ &\quad = \alpha_1 \mathbb{P}(B_n = 1 | (X_{n-1}, \ldots, X_1) \equiv 0) \\ &\qquad + \alpha_0 (1 - \mathbb{P}(B_n = 1 | (X_{n-1}, \ldots, X_1) \equiv 0)) \\ &\quad = \alpha_0 + (\alpha_1 - \alpha_0) \mathbb{P}(B_n = 1 | (X_{n-1}, \ldots, X_1) \equiv 0). \end{aligned}$$

Therefore

$$\begin{aligned} &\mathbb{P}(B_n = 1 | (X_{n-1}, \ldots, X_1) \equiv 0) \\ &\qquad = \frac{\mathbb{P}(X_n = 1 | (X_{n-1}, \ldots, X_1) \equiv 0) - \alpha_0}{\alpha_1 - \alpha_0} = \frac{A_n - \alpha_0}{\alpha_1 - \alpha_0}. \end{aligned}$$

Furthermore, using the above we get

$$\begin{aligned} &\mathbb{P}(X_n = 1, X_{n-1} = 0 | (X_{n-2}, \ldots, X_1) \equiv 0) \\ &\quad = \mathbb{P}(X_n = 1, X_{n-1} = 0 | B_{n-1} = 1, (X_{n-2}, \ldots, X_1) \equiv 0) \\ &\qquad \times \mathbb{P}(B_{n-1} = 1 | (X_{n-2}, \ldots, X_1) \equiv 0) \\ &\qquad + \mathbb{P}(X_n = 1, X_{n-1} = 0 | B_{n-1} = 0, (X_{n-2}, \ldots, X_1) \equiv 0) \\ &\qquad \times \mathbb{P}(B_{n-1} = 0 | (X_{n-2}, \ldots, X_1) \equiv 0) \end{aligned}$$

$$
\begin{aligned}
&= \mathbb{P}(X_n = 1 | X_{n-1} = 0, B_{n-1} = 1) \\
&\quad \times \mathbb{P}(X_{n-1} = 0 | B_{n-1} = 1) \frac{A_{n-1} - \alpha_0}{\alpha_1 - \alpha_0} \\
&\quad + \mathbb{P}(X_n = 1 | X_{n-1} = 0, B_{n-1} = 0) \\
&\quad\quad \times \mathbb{P}(X_{n-1} = 0 | B_{n-1} = 0) \left(1 - \frac{A_{n-1} - \alpha_0}{\alpha_1 - \alpha_0}\right) \\
&= [\alpha_1(1 - \gamma(1-p)) + \alpha_0\gamma(1-p)](1-\alpha_1) \frac{A_{n-1} - \alpha_0}{\alpha_1 - \alpha_0} \\
&\quad + [\alpha_1\gamma p + \alpha_0(1-\gamma p)](1-\alpha_0)\left(1 - \frac{A_{n-1} - \alpha_0}{\alpha_1 - \alpha_0}\right) \\
&= \frac{A_{n-1} - \alpha_0}{\alpha_1 - \alpha_0}([\alpha_1(1 - \gamma(1-p)) + \alpha_0\gamma(1-p)](1-\alpha_1) \\
&\qquad\qquad - [\alpha_1\gamma p + \alpha_0(1-\gamma p)](1-\alpha_0)) \\
&\quad + [\alpha_1\gamma p + \alpha_0(1-\gamma p)](1-\alpha_0).
\end{aligned}
$$

Finally observing that

$$
\begin{aligned}
&\frac{1}{\alpha_1 - \alpha_0}[[\alpha_1(1 - \gamma(1-p) + \alpha_0\gamma(1-p))](1-\alpha_1) \\
&\qquad\qquad - [\alpha_1\gamma p + \alpha_0(1-\gamma p)](1-\alpha_0)] \\
&\quad = \frac{1}{\alpha_1 - \alpha_0}[[\alpha_1 - (\alpha_1 - \alpha_0)\gamma(1-p)](1-\alpha_1) \\
&\qquad\qquad - [(\alpha_1 - \alpha_0)\gamma p + \alpha_0](1-\alpha_0)] \\
&\quad = \frac{1}{\alpha_1 - \alpha_0}[\alpha_1(1-\alpha_1) - \alpha_0(1-\alpha_0)] - \gamma(1-p)(1-\alpha_1) - \gamma p(1-\alpha_0) \\
&\quad = (1 - \alpha_0 - \alpha_1) - \gamma(1-p)(1-\alpha_1) - \gamma p(1-\alpha_0) \\
&\quad = (1 - \alpha_0 - \alpha_1) - \gamma[1 - \alpha_0 - (1-p)(\alpha_1 - \alpha_0)],
\end{aligned}
$$

and that

$$
\begin{aligned}
&-\frac{\alpha_0}{\alpha_1 - \alpha_0}([\alpha_1(1 - \gamma(1-p)) + \alpha_0\gamma(1-p)](1-\alpha_1) \\
&\qquad\qquad - [\alpha_1\gamma p + \alpha_0(1-\gamma p)](1-\alpha_0)) \\
&\quad + [\alpha_1\gamma p + \alpha_0(1-\gamma p)](1-\alpha_0) \\
&\quad = \frac{\alpha_1}{\alpha_1 - \alpha_0}[\alpha_1\gamma p + \alpha_0(1-\gamma p)](1-\alpha_0) \\
&\quad\quad - \frac{\alpha_0}{\alpha_1 - \alpha_0}[\alpha_1(1 - \gamma(1-p)) + \alpha_0\gamma(1-p)](1-\alpha_1)
\end{aligned}
$$

$$
\begin{aligned}
&= \frac{1}{\alpha_1 - \alpha_0}[[\alpha_0 + (\alpha_1 - \alpha_0)\gamma p](1-\alpha_0)\alpha_1 \\
&\qquad - [\alpha_1 - (\alpha_1 - \alpha_0)\gamma(1-p)](1-\alpha_1)\alpha_0] \\
&= \frac{1}{\alpha_1 - \alpha_0}[\alpha_0(1-\alpha_0)\alpha_1 - \alpha_1(1-\alpha_1)\alpha_0] \\
&\quad + \gamma p(1-\alpha_0)\alpha_1 + \gamma(1-p)(1-\alpha_1)\alpha_0 \\
&= \alpha_0\alpha_1 + \gamma p(1-\alpha_0)\alpha_1 + \gamma(1-p)(1-\alpha_1)\alpha_0 \\
&= \alpha_0\alpha_1 + \gamma[(p - \alpha_0 + (1-p)\alpha_0)\alpha_1 + (1-p)\alpha_0 - (1-p)\alpha_0\alpha_1] \\
&= \alpha_0\alpha_1 + \gamma[\alpha_1(1-\alpha_0) - (1-p)(\alpha_1 - \alpha_0)],
\end{aligned}
$$

completes the proof. □

We are now ready to prove Theorem 1.3.

Proof of Theorem 1.3. From Proposition 1.2 we know that the limit $A = \lim_{n\to\infty} A_n$ exists, and therefore we can take the limit of both sides of equation (12) ($A_n$ is easily seen to be uniformly bounded away from 1) to conclude that

$$A = \lim_{n\to\infty} A_n = \lim_{n\to\infty} \frac{CA_{n-1} + D}{1 - A_{n-1}} = \frac{CA + D}{1 - A}.$$

This gives us that

$$A - A^2 = CA + D,$$

and therefore

$$A^2 + (C-1)A + D = 0,$$

solving this equation gives

$$
\begin{aligned}
A &= \tfrac{1}{2}(1 - C \pm \sqrt{(1-C)^2 - 4D}) \\
&= \tfrac{1}{2}(1 - C \pm \sqrt{(2\alpha_1 - (1-C))^2 - 4(D + \alpha_1^2 - \alpha_1(1-C))}).
\end{aligned}
$$

We will now proceed to rule out one of the solutions:

$$
\begin{aligned}
&D + \alpha_1^2 - \alpha_1(1-C) \\
&\quad = \alpha_1^2 - \alpha_1 + \alpha_0\alpha_1 + \gamma(\alpha_1(1-\alpha_0) - (1-p)(\alpha_1 - \alpha_0)) \\
&\qquad + \alpha_1((1 - \alpha_0 - \alpha_1) - \gamma(1 - \alpha_0 - (1-p)(\alpha_1 - \alpha_0))) \\
&\quad = \gamma(\alpha_1(1-\alpha_0) - (1-p)(\alpha_1 - \alpha_0)) - \gamma\alpha_1(1 - \alpha_0 - (1-p)(\alpha_1 - \alpha_0)) \\
&\quad = \gamma(-(1-p)(\alpha_1 - \alpha_0) + (1-p)\alpha_1(\alpha_1 - \alpha_0)) \\
&\quad = -\gamma(1-p)(\alpha_1 - \alpha_0)(1 - \alpha_1).
\end{aligned}
$$

Using that $\gamma(1-p)((\alpha_1-\alpha_0)(1-\alpha_1)) \geq 0$ we get

$$\begin{aligned} &\tfrac{1}{2}(1-C+\sqrt{(2\alpha_1-(1-C))^2-4(D+\alpha_1^2-\alpha_1(1-C))}) \\ &\quad \geq \tfrac{1}{2}(1-C+\sqrt{(2\alpha_1-(1-C))^2}) \geq \tfrac{1}{2}(1-C+(2\alpha_1-(1-C))) = \alpha_1. \end{aligned}$$

Obviously we cannot have

$$A \geq \alpha_1,$$

since already [for $\gamma, p \in (0,1)$]

$$A_2 = \mathbb{P}(X_2 = 1 | X_1 = 0) < \alpha_1$$

and $A \leq A_n$ for every $n$. We conclude that

$$A = \tfrac{1}{2}(1-C-\sqrt{(1-C)^2-4D}).$$

Using Proposition 1.2 we then conclude that $p_{\max,\mu} = A$. Finally, the result for $p_{\min,\mu}$ follows from an easy symmetry argument. $\square$

Observe that when $\alpha_0 = \alpha_1 = \alpha$, $\{X_n\}_{n=1}^{\infty}$ is i.i.d. and $p_{\max,\mu} = p_{\min,\mu} = \alpha$. Note that in this case $C = 1 - 2\alpha - \gamma(1-\alpha)$ and $D = \alpha^2 + \gamma\alpha(1-\alpha)$, and so

$$\begin{aligned} p_{\max,\mu} &= \tfrac{1}{2}(2\alpha+\gamma(1-\alpha)-\sqrt{(2\alpha+\gamma(1-\alpha))^2-4(\alpha^2+\gamma\alpha(1-\alpha))}) \\ &= \tfrac{1}{2}(2\alpha+\gamma(1-\alpha)-\sqrt{(\gamma(1-\alpha))^2}) = \alpha, \end{aligned}$$

as we should get. Similarly one can check that $p_{\min,\mu} = \alpha$.

Furthermore if we choose $\gamma = 1$, $\{X_n\}_{n=1}^{\infty}$ is again i.i.d. and we would expect to get that $p_{\max,\mu} = p_{\min,\mu} = \alpha_0 + p(\alpha_1 - \alpha_0)$. Again, this is easy to check. Finally, as $\gamma \to 0$ we get that $p_{\max,\mu} \to \alpha_0$ and $p_{\min,\mu} \to \alpha_1$. It is not hard to see why this is what we should expect.

## 4. Continuous time domination results.

In this section we prove Theorem 1.4.

For $T > 0$, let $D_{\mathbb{N}}[0,T]$ be the set of functions from $[0,T]$ to $\mathbb{N}$ that are right-continuous and have left limits. Let $D_{\mathbb{N}}[0,\infty)$ be defined in the same way, but with $[0,T]$ replaced by $[0,\infty)$. Let a function be called a count path if it is a nondecreasing function that takes integer values and has jumps of size 1. Define $D_c \subset D_{\mathbb{N}}[0,1]$ to be the set of count paths. $D_c$ is closed under the Skorokhod topology, see [4], page 137.

Let $\alpha_0, \alpha_1, \gamma > 0$ and let $m$ be such that $\alpha_{0,m} := \alpha_0/m$, $\alpha_{1,m} := \alpha_1/m$, $\gamma_m := \gamma/m \in (0,1)$. Consider the model in the last section with $\alpha_0, \alpha_1$ and $\gamma$ replaced by $\alpha_{0,m}, \alpha_{1,m}$ and $\gamma_m$, respectively ($p$ is not changed). Denote the corresponding processes by $\{(B_n^m, X_n^m)\}_{n=1}^\infty$ but consider only the truncated part $\{(B_n^m, X_n^m)\}_{n=1}^m$. As in the Introduction, let

$$X_n^{c,m} = \sum_{i=1}^n X_i^m \qquad \text{for } n \in \{1, \dots, m\}.$$

Define the continuous time version $\{(B_t^m, X_t^m)\}_{t \in [0,1]}$ by letting

$$(B_t^m, X_t^m) = (B_n^m, X_n^{c,m}) \qquad \text{for } t \in [n-1, n)/m \text{ and } n \in \{1, \dots, m\}, \tag{13}$$

and $(B_{t=1}^m, X_{t=1}^m) = (B_m^m, X_m^{c,m})$. According to Theorem 1.3, we can couple the $\{(B_n^m, X_n^m)\}_{n=1}^m$ process with an i.i.d. process $\{Y_n^m\}_{n=1}^m$ with density $p_{\max,\mu_m}$ (where $\mu_m$ denotes the distribution of $\{X_n^m\}_{n=1}^\infty$) such that

$$Y_n^m \le X_n^m \qquad \forall n \in \{1, \dots, m\}. \tag{14}$$

Here $p_{\max,\mu_m}$ is given by Theorem 1.3. Define $\{Y_n^{c,m}\}_{n=1}^\infty$ in the obvious way and the continuous time version $\{Y_t^m\}_{t \in [0,1]}$ by letting

$$Y_t^m = Y_n^{c,m} \qquad \text{for } t \in [n-1, n)/m,\ n \in \{1, \dots, m\}$$

and $Y_{t=1}^m = Y_m^{c,m}$. We get from equation (14) that

$$X_t^m - Y_t^m \text{ is nondecreasing in } t \qquad \forall n \in \{1, \dots, m\}. \tag{15}$$

We state the following lemma; the proof is an elementary exercise in convergence in the Skorokhod topology.

LEMMA 4.1. *The set* $\{(f, g) \in D_c \times D_c : f - g$ *is nondecreasing*$\}$ *is closed in the product Skorokhod topology.*

Consider now $\{(B_t, X_t)\}_{t \in [0,1]}$ defined in Section 1. Recall that the flip rate intensities corresponding to $\{(B_t, X_t)\}_{t \in [0,1]}$ are

(16)

| from | to | with intensity |
|---|---|---|
| $(0, k)$ | $(1, k)$ | $\gamma p$ |
| $(1, k)$ | $(0, k)$ | $\gamma(1-p)$ |
| $(0, k)$ | $(0, k+1)$ | $\alpha_0$ |
| $(1, k)$ | $(1, k+1)$ | $\alpha_1$ |

for any $k \geq 0$. Observe that for $\{(B_n^m, X_n^{c,m})\}_{n=1}^{\infty}$ we have the transition probabilities

$$
\begin{array}{lll}
\text{from} & \text{to} & \text{w.p.} \\
(0,k) & (1,k) & (\gamma p/m)(1-\alpha_1/m) \\
(1,k) & (0,k) & (\gamma(1-p)/m)(1-\alpha_0/m) \\
(0,k) & (1,k+1) & \gamma p \alpha_1/m^2 \\
(1,k) & (0,k+1) & \gamma(1-p)\alpha_0/m^2 \\
(0,k) & (0,k) & (1-\alpha_0/m)(1-\gamma p/m) \\
(1,k) & (1,k) & (1-\alpha_1/m)(1-\gamma(1-p)/m) \\
(0,k) & (0,k+1) & (\alpha_0/m)(1-\gamma p/m) \\
(1,k) & (1,k+1) & (\alpha_1/m)(1-\gamma(1-p)/m).
\end{array} \tag{17}
$$

Using the flip rate intensities of equations (16) and (17), it is a standard result to show the next lemma. Again we omit the proof. However, see, for instance, [6] for a survey on the convergence of Markov processes in general.

Lemma 4.2. *The sequence of processes $\{(B_t^m, X_t^m)\}_{t\in[0,1]}$ defined above and indexed by $m$, converges weakly to the Markov process $\{(B_t, X_t)\}_{t\in[0,1]}$.*

We are now ready to prove our main results of this section. We will start by proving the following lemma.

Lemma 4.3. *With the assumptions of Theorem 1.4, we have that*

$$\lambda_{\max,\mu}(\alpha_0, \alpha_1, \gamma, p) \geq \bar{\lambda}.$$

Proof. We will start by constructing the coupling on the finite time interval $[0,1]$ and then argue that we can extend it to infinite time.

Let $\{(B_t^m, X_t^m, Y_t^m)\}_{t\in[0,1]}$ be any sequence of processes indexed by $m$ where, as indicated by the notation, the marginals $\{(B_t^m, X_t^m)\}_{t\in[0,1]}$ and $\{Y_t^m\}_{t\in[0,1]}$ have the distribution of the processes defined at the beginning of this section. Furthermore assume that these marginals are coupled so that $X_t^m - Y_t^m$ is nondecreasing for every $m$. Obviously the marginal $\{Y_t^m\}_{t\in[0,1]}$ converges weakly to a Poisson process $\{Y_t\}_{t\in[0,1]}$ with intensity

$$
\begin{aligned}
\lim_{m\to\infty} & m p_{\max,\mu_m} \\
&= \lim_{m\to\infty} \frac{1}{2}\Big(\alpha_0 + \alpha_1 + \gamma\Big(1 - \frac{1}{m}\alpha_0 - (1-p)\frac{1}{m}(\alpha_1 - \alpha_0)\Big)\Big) \\
&\quad - \frac{1}{2}\bigg(\Big(\alpha_0 + \alpha_1 + \gamma\Big(1 - \frac{1}{m}\alpha_0 - (1-p)\frac{1}{m}(\alpha_1-\alpha_0)\Big)\Big)^2 \\
&\qquad - 4\Big(\alpha_0\alpha_1 + \gamma\Big(\alpha_1\Big(1 - \frac{1}{m}\alpha_0\Big) - (1-p)(\alpha_1 - \alpha_0)\Big)\Big)\bigg)^{1/2}
\end{aligned}
$$

$$= \frac{1}{2}(\alpha_0 + \alpha_1 + \gamma - \sqrt{(\alpha_0 + \alpha_1 + \gamma)^2 - 4(\alpha_0\alpha_1 + \gamma(\alpha_0 + p(\alpha_1 - \alpha_0)))})$$
$$= \frac{1}{2}(\alpha_0 + \alpha_1 + \gamma - \sqrt{(\alpha_1 - \alpha_0 - \gamma)^2 + 4\gamma(1-p)(\alpha_1 - \alpha_0)}) = \bar{\lambda}.$$

Lemma 4.2 shows that also the sequence $\{(B_t^m, X_t^m)\}_{t\in[0,1]}$ converges weakly. It can then be argued that the sequence $\{(B_t^m, X_t^m, Y_t^m)\}_{t\in[0,1]}$ is tight and so there exists a subsequence $\{\{(B_t^{m(k)}, X_t^{m(k)}, Y_t^{m(k)})\}_{t\in[0,1]}\}_{k=1}^{\infty}$ that converges weakly to some process $\{(\tilde{B}_t, \tilde{X}_t, \tilde{Y}_t)\}_{t\in[0,1]}$. Of course, the marginal distribution $\{(\tilde{B}_t, \tilde{X}_t)\}_{t\in[0,1]}$ must be equal to the distribution of $\{(B_t, X_t)\}_{t\in[0,1]}$, and the marginal distribution $\{\tilde{Y}_t\}_{t\in[0,1]}$ must be equal to the distribution of $\{Y_t\}_{t\in[0,1]}$. Furthermore using Lemma 4.1 we conclude that

$$X_t - Y_t \text{ is nondecreasing.} \tag{18}$$

It not hard to see that we can adapt the proof to work for any time interval $[0,T]$. It is then easy to construct the coupling on $D_{\mathbb{N}}[0,\infty)$. Hence we have established that

$$\lambda_{\max,\mu}(\alpha_0, \alpha_1, \gamma, p) \geq \bar{\lambda}. \qquad \square$$

Considering $\{(B_n^m, X_n^m)\}_{n=1}^{\infty}$, let for every $m, i \geq 1$ $A_i^m := \mathbb{P}(X_i^m = 1 | X_{i-1}^m = \cdots = X_1^m = 0)$, and let $A^m := p_{\max,\mu_m} = \lim_{i\to\infty} A_i^m$. In our next lemma, we will need that $Tm$ (where $T > 0$) is an integer, which will not always be the case. However, adjusting the proofs for this is trivial and we therefore leave it to the reader. The same comment applies for other results to follow.

LEMMA 4.4. *For any $T > 0$,*

$$\lim_{m\to\infty} mA_{Tm}^m$$
$$= \bar{\lambda} + (p\alpha_1 + (1-p)\alpha_0 - \bar{\lambda}) \frac{e^{-(\alpha_0+\alpha_1+\gamma)T}}{e^{-\bar{\lambda}T}\mathbb{P}(X_t = 0 \ \forall t \in [0,T])}.$$

PROOF. Let $C^m, D^m$ denote $C, D$ of Theorem 1.3 with parameters $\alpha_0/m$, $\alpha_1/m$, $\gamma/m$ and $p$. By Lemma 3.2, for any $n$,

$$A_n^m - A^m = \frac{C^m A_{n-1}^m + D^m}{1 - A_{n-1}^m} - \frac{C^m A^m + D^m}{1 - A^m}$$
$$= \frac{(C^m A_{n-1}^m + D^m)(1 - A^m) - (C^m A^m + D^m)(1 - A_{n-1}^m)}{(1 - A_{n-1}^m)(1 - A^m)}$$
$$= \frac{C^m (A_{n-1}^m - A^m) + D^m (A_{n-1}^m - A^m)}{(1 - A_{n-1}^m)(1 - A^m)}$$

$$= (A^m_{n-1} - A^m)\frac{C^m + D^m}{(1 - A^m_{n-1})(1 - A^m)}$$

$$= \cdots = (A^m_1 - A^m)\Big(\frac{C^m + D^m}{1 - A^m}\Big)^{n-1}\frac{1}{\prod_{k=1}^{n-1}(1 - A^m_k)}.$$

Furthermore

$$\begin{aligned}C^m + D^m &= (1 - \alpha_0/m - \alpha_1/m) - \gamma/m(1 - \alpha_0/m - (1 - p)(\alpha_1/m - \alpha_0/m)) \\ &\quad + \alpha_0\alpha_1/m^2 + \gamma/m(\alpha_1/m(1 - \alpha_0/m) - (1 - p)(\alpha_1/m - \alpha_0/m)) \\ &= 1 - \alpha_0/m - \alpha_1/m + \alpha_0\alpha_1/m^2 \\ &\quad - \gamma/m(1 - \alpha_0/m - \alpha_1/m(1 - \alpha_0/m)) \\ &= (1 - \alpha_0/m)(1 - \alpha_1/m)(1 - \gamma/m).\end{aligned}$$

Recall also that we in Lemma 4.3 proved that $\bar{\lambda} = \lim_{m\to\infty} mA^m$. We get that:

1. $\lim_{m\to\infty}(C^m + D^m)^{Tm-1} = e^{-(\alpha_0+\alpha_1+\gamma)T}$,
2. $\lim_{m\to\infty}(1 - A^m)^{Tm-1} = \lim_{m\to\infty} e^{(Tm-1)\log(1-A^m)} = \lim_{m\to\infty} e^{(Tm-1)(-A^m+\mathcal{O}((A^m)^2))} = e^{-\bar{\lambda}T}$,
3. $\lim_{m\to\infty}\prod_{k=1}^{Tm-1}(1 - A^m_k) = \lim_{m\to\infty}\mathbb{P}(X^m_t = 0 \ \ \forall t \in [0, T - 1/m]) = \mathbb{P}(X_t = 0 \ \forall t \in [0, T])$,
4. $mA^m_1 = m(p\alpha_1/m + (1 - p)\alpha_0/m) = p\alpha_1 + (1 - p)\alpha_0$.

Therefore

$$\lim_{m\to\infty} mA^m_{Tm}$$

$$= \lim_{m\to\infty} mA^m + (mA^m_1 - mA^m)\Big(\frac{C^m + D^m}{1 - A^m}\Big)^{Tm-1}\frac{1}{\prod_{k=1}^{Tm-1}(1 - A^m_k)} \tag{19}$$

$$= \bar{\lambda} + (p\alpha_1 + (1 - p)\alpha_0 - \bar{\lambda})\frac{e^{-(\alpha_0+\alpha_1+\gamma)T}}{e^{-\bar{\lambda}T}\mathbb{P}(X_t = 0 \ \forall t \in [0, T])},$$

as desired. □

Next we prove the upper bound in Theorem 1.5 of $\lambda^T_{\max,\mu}(\alpha_0, \alpha_1, \gamma, p)$.

LEMMA 4.5. *For every choice of* $\alpha_0, \alpha_1, \gamma, T > 0$, *with* $\alpha_0 \le \alpha_1$ *and* $p \in (0, 1)$ *we have that there exists a constant* $E > 0$, *depending on* $\alpha_1, \alpha_0, \gamma$ *and* $p$ *such that*

$$\lambda^T_{\max,\mu}(\alpha_0, \alpha_1, \gamma, p) \le \bar{\lambda} + \frac{1}{T}(p\alpha_1 + (1 - p)\alpha_0 - \bar{\lambda})\frac{1 - e^{-TE}}{E}.$$

PROOF. We have that

$$\mathbb{P}(X_t^m = 0\ \forall t \in [0,T]) = \prod_{k=1}^{Tm}(1 - A_k^m) = e^{\sum_{k=1}^{Tm}\log(1-A_k^m)}$$

$$(20) \qquad = e^{-\sum_{k=1}^{Tm} A_k^m + \mathcal{O}((A_k^m)^2)} = e^{\mathcal{O}(1/m) - \sum_{k=1}^{Tm} A_k^m}.$$

Using equation (20) it is easy to see that it suffices to get an estimate on $\sum_{k=1}^{Tm} A_k^m$. To that end, let $n > 0$ be an integer and let $T_k := kT/n$ for $k \in \{1, \ldots, n\}$. Using that for fixed $m$, $A_k^m$ is decreasing in $k$, we get that

$$\sum_{k=1}^{Tm} A_k^m = \sum_{k=1}^{T_1 m} A_k^m + \sum_{k=T_1 m+1}^{T_2 m} A_k^m + \cdots + \sum_{k=T_{n-1}m+1}^{T_n m} A_k^m$$

$$\le T_1 m A_1^m + (T_2 - T_1) m A_{T_1 m}^m + \cdots + (T_n - T_{n-1}) m A_{T_{n-1}m}^m.$$

Using equation (19), that $mA_1^m = p\alpha_1 + (1 - p\alpha_0)$ and that $(T_k - T_{k-1}) = T/n$ for every $k$, we get that

$$\lim_{m\to\infty} \sum_{k=1}^{Tm} A_k^m$$

$$(21) \qquad \le \frac{T}{n}(p\alpha_1 + (1-p)\alpha_0)$$

$$+ \sum_{k=1}^{n-1} \frac{T}{n}\left[\bar{\lambda} + (p\alpha_1 + (1-p)\alpha_0 - \bar{\lambda})\frac{e^{-(\alpha_0+\alpha_1+\gamma)T_k}}{e^{-\bar{\lambda}T_k}\mathbb{P}(X_t = 0\ \forall t \in [0,T_k])}\right].$$

Note that the existence of this limit follows from the existence of the limit on the left-hand side of equation (20). We observe that trivially $\mathbb{P}(X_t = 0\ \forall t \in [0,T]) \ge e^{-\alpha_1 T}$ and so we get that

$$\frac{e^{-(\alpha_0+\alpha_1+\gamma)T_k}}{e^{-\bar{\lambda}T_k}\mathbb{P}(X_t = 0\ \forall t \in [0,T_k])}$$

$$\le \exp(-(\alpha_0 + \alpha_1 + \gamma)T_k + \bar{\lambda}T_k + \alpha_1 T_k)$$

$$= \exp\biggl(\frac{-T_k}{2}(\alpha_0 + \alpha_1 + \gamma$$

$$+ \sqrt{(\alpha_1 - \alpha_0 - \gamma)^2 + 4\gamma(1-p)(\alpha_1 - \alpha_0)}) + \alpha_1 T_k\biggr)$$

$$= \exp\left(\frac{-T_k}{2}(\alpha_0 + \alpha_1 + \gamma + |\alpha_1 - \alpha_0 - \gamma| + 2E) + \alpha_1 T_k\right)$$

$$= \exp(T_k(\alpha_1 - \max(\alpha_1, \alpha_0 + \gamma) - E)) \le e^{-ET_k},$$

where $E$ solves the equation

$$|\alpha_1 - \alpha_0 - \gamma| + 2E = \sqrt{|\alpha_1 - \alpha_0 - \gamma|^2 + 4\gamma(1-p)(\alpha_1 - \alpha_0)}.$$

We get that

$$\begin{aligned}
&\sum_{k=1}^{n-1} \frac{T}{n} \frac{e^{-(\alpha_0+\alpha_1+\gamma)T_k}}{e^{-\bar{\lambda}T_k}\mathbb{P}(X_t = 0\ \forall t \in [0,T_k])} \\
&\quad \le \frac{T}{n}\sum_{k=1}^{n-1} e^{-ET_k} = \frac{T}{n}\sum_{k=1}^{n-1}(e^{-ET/n})^k \qquad (22)\\
&\quad = \frac{T}{n}\left(\frac{1-e^{-TE}}{1-e^{-TE/n}} - 1\right) = \frac{T}{n}\left(\frac{e^{-TE/n} - e^{-TE}}{1-e^{-TE/n}}\right) \\
&\quad = \frac{T}{n}\left(\frac{e^{-TE/n} - e^{-TE}}{TE/n + \mathcal{O}(1/n^2)}\right) = \left(\frac{e^{-TE/n} - e^{-TE}}{E + \mathcal{O}(1/n)}\right).
\end{aligned}$$

Combining equations (21) and (22) and taking the limit as $n$ tends to infinity (after taking the limit as $m$ tends to infinity), we get that

$$\lim_{m\to\infty} \sum_{k=1}^{Tm} A_k^m \le T\bar{\lambda} + (p\alpha_1 + (1-p)\alpha_0 - \bar{\lambda})\frac{1-e^{-TE/2}}{E}.$$

Combining equation (20) with this yields

$$\begin{aligned}
&\mathbb{P}(X_t = 0\ \forall t \in [0,T]) \\
&\quad = \lim_{m\to\infty} \mathbb{P}(X_t^m = 0\ \forall t \in [0,T]) \\
&\quad \ge \exp\left(-\left(T\bar{\lambda} + (p\alpha_1 + (1-p)\alpha_0 - \bar{\lambda})\frac{1-e^{-TE/2}}{E}\right)\right).
\end{aligned}$$

Finally we conclude that

$$\lambda^T_{\max,\mu}(\alpha_0,\alpha_1,\gamma,p) \le \bar{\lambda} + \frac{1}{T}(p\alpha_1 + (1-p)\alpha_0 - \bar{\lambda})\frac{1-e^{-TE/2}}{E}. \qquad \square$$

REMARK. It is interesting that in the above proof we "lift ourselves up by the boots" by using a simple estimate for $\mathbb{P}(X_t = 0\ \forall t \in [0,T])$ to obtain a better one.

PROOF OF THEOREM 1.4. The first statement follows immediately from Lemmas 4.3 and 4.5, by letting $T$ tend to infinity.

We can of course trivially conclude that $\lambda_{\min,\mu}(\alpha_0,\alpha_1,\gamma,p) \le \alpha_1$. To see why we have equality consider the event

$$\{\text{There are at least } k \text{ arrivals during } [0,1]\}.$$

Let $\alpha < \alpha_1$, we have that

$$\mathrm{Poi}_\alpha(\text{There are at least } k \text{ arrivals during } [0,1]) = \sum_{l=k}^{\infty} e^{-\alpha}\frac{\alpha^l}{l!}.$$

We also see that

$$\mathrm{Poi}^{\gamma,p}_{\alpha_0,\alpha_1}(\text{There are at least } k \text{ arrivals during } [0,1]) \geq pe^{-\gamma}\sum_{l=k}^{\infty} e^{-\alpha_1}\frac{\alpha_1^l}{l!}.$$

Since

$$\frac{\sum_{l=k}^{\infty} e^{-\alpha_1}\alpha_1^l/l!}{\sum_{l=k}^{\infty} e^{-\alpha}\alpha^l/l!} \longrightarrow_{k\to\infty} \infty,$$

we get that for every $\alpha < \alpha_1$, $\gamma > 0$ and $p > 0$ there exists a $k$ such that

$$\mathrm{Poi}_\alpha(\text{There are at least } k \text{ arrivals during } [0,1])$$
$$< \mathrm{Poi}^{\gamma,p}_{\alpha_0,\alpha_1}(\text{There are at least } k \text{ arrivals during } [0,1]).$$

Obviously this contradicts

$$\mathrm{Poi}^{\gamma,p}_{\alpha_0,\alpha_1} \preceq \mathrm{Poi}_\alpha,$$

and so $\lambda_{\min,\mu}(\alpha_0,\alpha_1,\gamma,p) \geq \alpha_1$. □

REMARK. It is actually possible to prove the statement without using Lemma 4.5, and we give here a short informal description how this can be done. To show that $\lambda_{\min,\mu} \leq \bar{\lambda}$ directly from Theorem 1.3, start with the processes $\{(B_t,X_t)\}_{t\geq 0}$ and $\{Y_t\}_{t\geq 0}$ with distributions as indicated by the notation (the latter process with parameter $\lambda$) coupled so that $\{X_t\}_{t\geq 0}$ has an arrival whenever $\{Y_t\}_{t\geq 0}$ has an arrival. For any $m$, it is straightforward to discretize these processes resulting in processes $\{(B_n^m,X_n^m)\}_{n=1}^\infty$ and $\{Y_n^m\}_{n=1}^\infty$ with distributions as indicated by the notation, but with parameters $\alpha_0/m + \mathcal{O}(1/m^2), \alpha_1/m + \mathcal{O}(1/m^2), \gamma/m + \mathcal{O}(1/m^2), p$ and $\lambda/m + \mathcal{O}(1/m^2)$, respectively. Furthermore this is done so that $\{X_n^m\}_{n=1}^\infty$ is coupled above $\{Y_n^m\}_{n=1}^\infty$. Using Theorem 1.3 we arrive at

$$\lambda/m + \mathcal{O}(1/m^2) \leq p_{\max,\mu^m},$$

where $\mu^m$ is the distribution of $\{X_n^m\}_{n=1}^\infty$. Multiplying with $m$ and letting $m$ go to infinity gives the result.

In the next section we will need the following easy corollary to Theorem 1.4.

COROLLARY 4.6. *For any $\delta < \min(\delta_1, \delta_0+\gamma)$ we can find a $0<p<1$ close enough to one so that*

$$\mathrm{Poi}_\delta \preceq \mathrm{Poi}^{\gamma,p}_{\delta_0,\delta_1}.$$

PROOF. We just need to observe that

$$
\begin{aligned}
&\lim_{p\to 1} \lambda_{\max,\mu}(\delta_0,\delta_1,\gamma,p)\\
&\qquad = \lim_{p\to 1} \tfrac{1}{2}\Bigl(\delta_0+\delta_1+\gamma-\sqrt{(\delta_1-\delta_0-\gamma)^2+4\gamma(1-p)(\delta_1-\delta_0)}\Bigr) \qquad (23)\\
&\qquad = \tfrac{1}{2}(\delta_1+\delta_0+\gamma-|\delta_1-\delta_0-\gamma|) = \min(\delta_1,\delta_0+\gamma). \qquad \square
\end{aligned}
$$

Observe that

$$\mathrm{Poi}^{\gamma,p}_{\delta_0,\delta_1}(\text{There are no arrivals in } [0,t]) \ge (1-p)e^{-\gamma t}e^{-\delta_0 t} = (1-p)e^{-(\gamma+\delta_0)t}$$

and that

$$\mathrm{Poi}_\delta(\text{There are no arrivals in } [0,t]) = e^{-\delta t}.$$

Therefore, if $\delta > \gamma+\delta_0$, we have for fixed $p$ and some $t$ that

$$e^{-\delta t} \le (1-p)e^{-(\gamma+\delta_0)t},$$

and so we cannot have that

$$\mathrm{Poi}_\delta \preceq \mathrm{Poi}^{\gamma,p}_{\delta_0,\delta_1},$$

which is an alternative way to see why the limit in equation (23) cannot simply be equal to $\delta_1$.

PROOF OF THEOREM 1.5. The upper bound is just Lemma 4.5.

For the lower bound, we start by observing that using Theorem 1.4 we trivially get that $\mathbb{P}(X_t = 0\ \forall t \in [0,T]) \le e^{-\bar\lambda T}$. Therefore by equation (19), using that $p\alpha_1 + (1-p)\alpha_0 \ge \bar\lambda$ [which follows from the fact that $\bar\lambda$ is increasing in $\gamma$ with limit $p\alpha_1+(1-p)\alpha_0$, see remark after Theorem 1.4],

$$
\begin{aligned}
&\lim_{m\to\infty} mA^m_{Tm}\\
&\qquad \ge \bar\lambda + (p\alpha_1+(1-p)\alpha_0-\bar\lambda)\frac{e^{-(\alpha_0+\alpha_1+\gamma)T}}{e^{-2T\bar\lambda}}\\
&\qquad = \bar\lambda + (p\alpha_1+(1-p)\alpha_0-\bar\lambda)e^{-TL}.
\end{aligned}
$$

We therefore need to show that $\lambda^T_{\max,\mu} \ge \lim_{m\to\infty} mA^m_{Tm}$. Observe that the second marginal of the discrete time process $\{(B^m_n, X^m_n)\}_{n=1}^{Tm}$ trivially dominates an i.i.d. sequence of density $A^m_{Tm}$. Therefore, going through a limiting procedure very similar to the one of Lemma 4.3, we get that the second marginal of $\{(B_t,X_t)\}_{t\in[0,T]}$ dominates a Poisson process with intensity $\lim_{m\to\infty} mA^m_{Tm}$ on the time interval $[0,T]$.

The result for $\lambda^T_{\min,\mu}$ follows as in the proof of Theorem 1.4. $\square$

## 5. CPREE-results.

PROOF OF THEOREM 1.6. For every site $s \in S$ the recoveries of the $\Psi^{\gamma,p,A}_{\delta_0,\delta_1}$ process at that site has the same distribution as the arrivals of a $\mathrm{Poi}^{\gamma,p}_{\delta_0,\delta_1}$ process. By Corollary 4.6 we can couple the processes $\Psi^{\gamma,p,A}_{\delta_0,\delta_1}$ and $\Psi^{A}_{\delta}$ so that at every site the former has a recovery whenever the latter does. Furthermore, coupling the infection rates are of course trivial. This gives the first result [equation (11)].

Trivially the statement cannot hold if $\delta > \delta_1$. Furthermore, if $\delta > \delta_0 + \gamma$, then by letting $x \in A$, and noting that for $T$ large enough,

$$\begin{aligned}&\Psi^{\gamma,p,A}_{\delta_0,\delta_1}(\sigma_t(x) = 1\ \forall t \in [0,T])\\ &\qquad \geq (1-p)e^{-\delta_0 T}e^{-\gamma T} > e^{-\delta T} \geq \Psi^{A}_{\delta}(\sigma_t(x) = 1\ \forall t \in [0,T]),\end{aligned}$$

we are done. □

For $A \subset S$ such that $|A| < \infty$, let $\Psi^{\gamma,p,A}_{\delta_0,\infty,B(A)\equiv 0}$ denote the CPREE where a site $s \in S$ *always* is healthy (i.e., in state 0) as long as the background process of the site $s$ is in state 1. That is, we do not allow the site to become infected if the background process of the site is in state 1. More precisely, for any graph $G = (S,E)$ let $\{(B_t, Y_t)\}_{t\geq 0}$ be a pair of processes with state space $\{\{0,1\}\times\{0,1\}\}^S$ such that $B_0 \sim \pi_p$ conditioned on the event that $B_0(s) = 0$ for every $s \in A$, and let $Y_0(s) = 1$ iff $s \in A$. Observe that the conditioning does not affect the probability of $B_0(s)$ being 0 or 1 for any $s \notin A$. Let the pair evolve according to the following flip rate intensities at any site $s$.

| from | to | with intensity |
|---|---|---|
| $(0,0)$ | $(1,0)$ | $\gamma p$ |
| $(0,1)$ | $(1,0)$ | $\gamma p$ |
| $(1,0)$ | $(0,0)$ | $\gamma(1-p)$ |
| $(0,0)$ | $(0,1)$ | $\sum_{(s',s)\in E} Y_t(s')$ |
| $(0,1)$ | $(0,0)$ | $\delta_0$. |

Observe that with this definition the state $(1,1)$ is not allowed. Informally, this can be interpreted as letting the rate of recovery when $B_t(s) = 1$ be infinite, hence the notation.

PROOF OF THEOREM 1.7. We start by observing that it is easy to see from the definitions of $\Psi^{\gamma,p,A}_{\delta_0,\delta_1,B_0(A)\equiv 0}$ and $\Psi^{\gamma,p,A}_{\delta_0,\infty,B(A)\equiv 0}$ that

$$\Psi^{\gamma,p,A}_{\delta_0,\infty,B(A)\equiv 0} \preceq \Psi^{\gamma,p,A}_{\delta_0,\delta_1,B_0(A)\equiv 0}.$$

We will construct $\{(B_t, Y_t)\}_{t\geq 0}$ to have distribution $\Psi^{\gamma,p,A}_{\delta_0,\infty,B(A)\equiv 0}$ for some $p$ close to 0, and couple it with a process $\{Y'_t\}_{t\geq 0}$ such that $\{Y_t\}_{t\geq 0}$ stochastically dominates $\{Y'_t\}_{t\geq 0}$. It will be easy to see that in turn $\{Y'_t\}_{t\geq 0}$ will stochastically dominate $\Psi^{\lambda,A}_{\delta}$.

We now proceed to the actual construction. Let $B_0 \sim \pi_p$, conditioned on the event that $B_0(s) = 0$ for every $s \in A$. For every site $s \in S$, associate an independent process $\{B_t(s), X_t(s)\}_{t\geq 0}$ such that $\{1 - B_t(s), X_t(s)\}_{t\geq 0}$ is the model of Theorem 1.4 with $\alpha_0 = 0$, $\alpha_1 = \Delta_G$ and with $p$ replaced by $1 - p$. We get from Theorem 1.4 that

$$\begin{aligned}
&\lambda_{\max,\mu}(0, \Delta_G, \gamma, 1-p)\\
&\qquad = \tfrac{1}{2}\Big(0 + \Delta_G + \gamma - \sqrt{(\Delta_G - 0 - \gamma)^2 + 4\gamma(1-(1-p))(\Delta_G - 0)}\Big)\\
&\qquad = \tfrac{1}{2}\Big(\Delta_G + \gamma - \sqrt{(\Delta_G - \gamma)^2 + 4\Delta_G\gamma p}\Big),
\end{aligned}$$

that is, we can couple the pair of processes $\{B_t(s), X_t(s)\}_{t\geq 0}$ with a Poisson process $\{X'_t(s)\}_{t\geq 0}$ with intensity $\lambda_{\max,\mu}(0, \Delta_G, \gamma, 1-p)$ such that if this latter process has an arrival then so does $\{X_t(s)\}_{t\geq 0}$. There is a slight issue with $s \in A$, where we have conditioned that $B_0(s) = 0$. However, this corresponds to conditioning that the background process of Theorem 1.4 starts in state 1, and it is not hard to see that the conclusion of the theorem is still valid in this case. Informally, if we in this theorem start with the background process in state 1, this means that we are starting in the higher intensity state, and so it becomes "easier to dominate." It is easy to make this statement precise.

Let for every $s \in S$, $\{D_t(s)\}_{t\geq 0}$ be a Poisson process with intensity $\delta_0$, independent of each other and all other processes, and consider some quadruple $\{B_t(s), X_t(s), X'_t(s), D_t(s)\}_{t\geq 0}$ with marginal distributions as indicated by the notation. We now proceed to construct $\{(B_t, Y_t)\}_{t\geq 0}$ (the first marginal is of course already defined) and $\{Y'_t\}_{t\geq 0}$ from these four processes. Let $Y_0(s) = Y'_0(s)$ for every $s \in S$ and let $Y_0(s) = Y'_0(s) = 1$ iff $s \in A$. Let for every $s \in S$ $\{(B_t(s), Y_t(s))\}_{t\geq 0}$ and $\{Y'_t(s)\}_{t\geq 0}$ denote the marginals of the processes $\{(B_t, Y_t)\}_{t\geq 0}$ and $\{Y'_t\}_{t\geq 0}$ at the site $s$.

Let $N(s, Y_\tau(s)), N(s, Y'_\tau(s))$ denote the number of neighbors of the site $s$ that are infected at time $\tau$ under $\{Y_t\}_{t\geq 0}$ and $\{Y'_t\}_{t\geq 0}$, respectively. Recall that by definition, for any $s \in S$, $Y_0(s) = Y'_0(s) = 0$ if $B_0(s) = 1$. We will write $X_\tau(s) \neq X_{\tau^-}(s)$, $X'_\tau(s) \neq X'_{\tau^-}(s)$ and $D_\tau(s) \neq D_{\tau^-}(s)$ to indicate that these processes have an arrival at time $\tau$. Observe that by construction, for every $s \in S$ and $t \geq 0$, if $X'_\tau(s) \neq X'_{\tau^-}(s)$ then $X_\tau(s) \neq X_{\tau^-}(s)$. We will also write $B_{\tau^-}(s) < B_\tau(s)$ when we mean that the $B_t$ process flips from 0 to 1 at time $\tau$.

At time $\tau$, $\{(Y_t(s), Y'_t(s))\}_{t\geq 0}$ will change:

$$
\text{(24)}\qquad \begin{array}{lll} \text{from} & \text{to} & \text{if} \\ (1,1) & (0,0) & D_\tau(s)\neq D_{\tau^-}(s) \text{ or } B_{\tau^-}(s)<B_\tau(s) \\ (1,0) & (0,0) & D_\tau(s)\neq D_{\tau^-}(s) \text{ or } B_{\tau^-}(s)<B_\tau(s) \end{array}
$$

and also:

$$
\text{(25)}\qquad \begin{array}{llcl} \text{from} & \text{to} & \text{w.p.} & \text{if} \\ (0,0) & (1,1) & N(s,Y'_\tau(s))/\Delta_G & X'_\tau(s)\neq X'_{\tau^-}(s) \\ (0,0) & (1,0) & (N(s,Y_\tau(s))-N(s,Y'_\tau(s)))/\Delta_G & X'_\tau(s)\neq X'_{\tau^-}(s) \\ (0,0) & (1,0) & N(s,Y_\tau(s))/\Delta_G & X'_\tau(s)= X'_{\tau^-}(s) \\ & & & X_\tau(s)\neq X_{\tau^-}(s) \\ (1,0) & (1,1) & N(s,Y'_\tau(s))/\Delta_G & X'_\tau(s)\neq X'_{\tau^-}(s). \end{array}
$$

No other transitions are allowed. Note that by construction $\{X_t(s)\}_{t\geq 0}$ and $\{X'_t(s)\}_{t\geq 0}$ only have arrivals when $\{B_t(s)\}_{t\geq 0}$ is in state 0. Therefore, these rates make sure that $\{Y_t\}_{t\geq 0}$ and $\{Y'_t\}_{t\geq 0}$ are in state 0 when $\{B_t(s)\}_{t\geq 0}$ is in state 1. Note also that since $N(s, Y'_0(s)) = N(s, Y_0(s))$ for every $s\in S$, the rates make sure that

$$
Y'_t(s)\leq Y_t(s) \qquad \forall s\in S,\ t\geq 0,
$$

and that $N(s, Y'_t(s))\leq N(s, Y_t(s))$ for every $s\in S$ and $t\geq 0$.

It remains to check that $\{(B_t, Y_t)\}_{t\geq 0}$ and $\{Y'_t\}_{t\geq 0}$ have the right distribution. As noted above $\{Y_t\}_{t\geq 0}$ is 0 if $\{B_t\}_{t\geq 0}$ is 1. Furthermore it is easy to see that when $\{B_t\}_{t\geq 0}$ is 0, $\{Y_t\}_{t\geq 0}$ flips from 0 to 1 at rate $N(s, Y_\tau(s))$ and from 1 to 0 at rate $\delta_0$. It is also easy to see that $\{Y'_t\}_{t\geq 0}$ flips from 1 to 0 at a rate which is the minimum of two exponentially distributed times with parameters $\delta_0$ and $\gamma p$, the latter being the rate at which $\{B_t\}_{t\geq 0}$ flips from 0 to 1. Hence $\{Y'_t\}_{t\geq 0}$ flips from 1 to 0 at rate $\delta_0+\gamma p$ and by choosing p small enough this is less that $\delta$. It also not hard to see that $\{Y'_t\}_{t\geq 0}$ flips from 0 to 1 at a rate $\lambda_{\max,\mu}(0,\Delta_G,\gamma,1-p)N(s,Y'_t(s))/\Delta_G$. Furthermore by choosing $p$ perhaps even smaller, we get that

$$
\begin{aligned}
&\lambda_{\max,\mu}(0,\Delta_G,\gamma,1-p)N(s,Y'_t(s))/\Delta_G \\
&\quad = \frac{N(s,Y'_t(s))}{2\Delta_G}\Bigl(\Delta_G+\gamma-\sqrt{(\Delta_G-\gamma)^2+4\Delta_G\gamma p}\Bigr) \\
&\quad \geq \lambda N(s,Y'_t(s)).
\end{aligned}
$$

Here we used that $\gamma\geq\Delta_G$. Therefore $\{Y'_t\}_{t\geq 0}$ is a contact process with infection rate larger than $\lambda$ and with recovery rate less that $\delta$, and so the distribution of $\{Y'_t\}_{t\geq 0}$ dominates $\Psi_\delta^{\lambda,A}$. $\square$

EXAMPLE 5.1. Let $S = \mathbb{Z}$ and $A = \{-n, \ldots, 0\}$. We have that

$$\Psi^{\gamma,p,A}_{\delta_0,\infty,B_0(A)\equiv 0}(\sigma_t(1) = 0\ \forall t \in [0,T]) \\ \geq \Psi^{\gamma,p,A}_{\delta_0,\infty,B_0(A)\equiv 0}(B_t(1) = 1\ \forall t \in [0,T]) = pe^{-\gamma T}.$$

Furthermore, using Corollary 3.22 and Theorem 3.29 of [10], it follows after some work that for some constants $K, \varepsilon, \varepsilon' > 0$,

$$\Psi^{\lambda,A}_{\delta}(\sigma_t(1) = 0\ \forall t \in [0,T]) \leq Ke^{-\varepsilon n} + e^{-\varepsilon' T}.$$

By letting $n$ go to infinity we see that if $\gamma < \varepsilon'$ we cannot have

$$\Psi^{\lambda,A}_{\delta} \preceq \Psi^{\gamma,p,A}_{\delta_0,\infty,B_0(A)\equiv 0}.$$

Furthermore it is possible to modify this example to work for large $\delta_1$ rather than $\delta_1 = \infty$. This is done by considering how long the site $\{1\}$ is infected during the time interval $[0,T]$ rather than the probability of this site not being infected at all during $[0,T]$.

EXAMPLE 5.2. Assume that $x, y \in S$, are neighbors and that $\lambda = 1$. We have that $\Psi^{\lambda,\{x\}}_{\delta}(\sigma_t(y) = 1) = t + o(t)$ while $\Psi^{\gamma,p,\{x\}}_{\delta_0,\infty,B_0(x)=0}(\sigma_t(y) = 1) = (1-p)t + o(t)$, hence for $p$ positive we can find $t$ small enough so that $\Psi^{\lambda,\{x\}}_{\delta}(\sigma_t(y) = 1) > \Psi^{\gamma,p,\{x\}}_{\delta_0,\infty,B_0(x)=0}(\sigma_t(y) = 1)$.

We are now ready to prove Theorem 1.8.

PROOF OF THEOREM 1.8. We will start with the existence of $p_{c1}$ and $p_{c2}$.

Let $0 < p_1 \leq p_2 < \infty$ and let $\{B^1_t\}_{t\geq 0}$, $\{B^2_t\}_{t\geq 0}$ be two background processes with parameters $p_1, p_2$, respectively. Let $B^1_0$ have distribution $\pi_{p_1}$ and $B^2_0$ have distribution $\pi_{p_2}$ and couple them so that $B^1_0(s) \leq B^2_0(s)$ for every $s \in S$. It is easy to see that we can then couple the processes so that

$$B^1_t(s) \leq B^2_t(s) \qquad \forall t \geq 0,\ \forall s \in S.$$

Using these processes to construct $\{(B^1_t, Y^1_t)\}_{t\geq 0}$ and $\{(B^2_t, Y^2_t)\}_{t\geq 0}$ with distributions $\Psi^{\gamma,p_1,\{s\}}_{\delta_0,\delta_1}$ and $\Psi^{\gamma,p_2,\{s\}}_{\delta_0,\delta_1}$, it is easy to see that we can couple the marginals $\{Y^1_t\}_{t\geq 0}$, $\{Y^2_t\}_{t\geq 0}$ so that

$$Y^2_t(s) \leq Y^1_t(s) \qquad \forall t \geq 0,\ \forall s \in S.$$

This establishes the existence of $p_{c1}$ and $p_{c2}$.

Consider now the part of statement 4 concerning $p_{c2} > 0$. Choose $\delta > \delta_0$ close enough to $\delta_0$ and $\lambda < 1$ close enough to 1 so that the contact process $\Psi_\delta^{\gamma,\{s\}}$ survives weakly. Observing that $\Psi_{\delta_0,\infty}^{\gamma,p,\{s\}}$ is a convex combination of $\Psi_{\delta_0,\infty,B(s)=0}^{\gamma,p,\{s\}}$ and $\Psi_{\delta_0,\infty,B(s)=1}^{\gamma,p,A}$ and using Theorem 1.7 gives the result.

All of the statements about $p_{c1}, p_{c2} > 0$ are proved in exactly the same way. All of the statements about $p_{c1}, p_{c2} < 1$ are proved in a similar way, but follow even easier since we can use Theorem 1.6 directly without worrying about the initial state of the background process at $s$. $\square$

Proof of Proposition 1.9. We will show the theorem for $p_{c2}$, the proof for $p_{c1}$ is identical. First, we use Taylor's expansion to see that

$$\begin{aligned}
&\lim_{\gamma\to\infty} \lambda_{\max,\mu}(\delta_0,\delta_1,\gamma,p)\\
&\quad= \lim_{\gamma\to\infty} \frac{1}{2}\left(\delta_0+\delta_1+\gamma-\sqrt{(\delta_0+\delta_1+\gamma)^2-4(\delta_0\delta_1+\gamma(\delta_0+p(\delta_1-\delta_0)))}\right)\\
&\quad= \lim_{\gamma\to\infty} \frac{1}{2}(\delta_0+\delta_1+\gamma)\left(1-\sqrt{1-4\frac{\delta_0\delta_1+\gamma(\delta_0+p(\delta_1-\delta_0))}{(\delta_0+\delta_1+\gamma)^2}}\right)\\
&\quad= \lim_{\gamma\to\infty} \frac{1}{2}(\delta_0+\delta_1+\gamma)\\
&\qquad\times\left(1-\left(1-\frac{4(\delta_0\delta_1+\gamma(\delta_0+p(\delta_1-\delta_0)))/(\delta_0+\delta_1+\gamma)^2}{2}\right)\right.\\
&\qquad\qquad\left.+\mathcal{O}\left(\frac{1}{\gamma^2}\right)\right)\\
&\quad= \lim_{\gamma\to\infty}\frac{\delta_0\delta_1+\gamma(\delta_0+p(\delta_1-\delta_0))}{\delta_0+\delta_1+\gamma}+\mathcal{O}\left(\frac{1}{\gamma}\right)=\delta_0+p(\delta_1-\delta_0).
\end{aligned}$$

It is now clear from Theorem 1.4, the proof of Theorem 1.6 and the above calculation that given any $\varepsilon > 0$, we can find $\gamma'$ large enough so that with $\delta = \delta_0 + p(\delta_1 - \delta_0) - \varepsilon$ we have that for all $\gamma \geq \gamma'$

$$\Psi_{\delta_0,\delta_1}^{\gamma,p,A} \preceq \Psi_\delta^A$$

and so the $\Psi_{\delta_0,\delta_1}^{\gamma,p,A}$ dies out if $\delta_0 + p(\delta_1 - \delta_0) - \varepsilon > \delta_{c2}$. This is the same as saying that for any $\varepsilon > 0$ there exists $\gamma'$ large enough so that for all $\gamma \geq \gamma'$, if

$$p > \frac{\delta_{c2} - \delta_0 + \varepsilon}{\delta_1 - \delta_0},$$

the process dies out. Therefore for every $\gamma \geq \gamma'$ we have that

$$p_{c2}(\delta_0, \delta_1, \gamma) \leq \frac{\delta_{c2} - \delta_0 + \varepsilon}{\delta_1 - \delta_0}.$$

We can therefore conclude that

$$\limsup_{\gamma\to\infty} p_{c2}(\delta_0,\delta_1,\gamma) \le \frac{\delta_{c2}-\delta_0}{\delta_1-\delta_0}. \qquad \square$$

PROOF OF PROPOSITION 1.10. We show the proposition for $p_{c2}$; the proof for $p_{c1}$ is identical. Using the trivial facts that $\sqrt{1-x} \le 1 - x/2$ for $0 \le x \le 1$ and that

$$0 \le \frac{4p\gamma(\delta_1-\delta_0)}{(\delta_0-\delta_1-\gamma)^2} \le 1,$$

we get that for any $p$,

$$\begin{aligned}
&\lambda_{\max,\mu}(\delta_0,\delta_1,\gamma,p)\\
&\quad = \frac{1}{2}(\delta_0+\delta_1+\gamma-\sqrt{(\delta_1-\delta_0-\gamma)^2+4\gamma(1-p)(\delta_1-\delta_0)})\\
&\quad = \frac{1}{2}(\delta_0+\delta_1+\gamma-\sqrt{(\delta_0-\delta_1-\gamma)^2-4p\gamma(\delta_1-\delta_0)})\\
&\quad = \frac{1}{2}\left(\delta_0+\delta_1+\gamma-|\delta_0-\delta_1-\gamma|\sqrt{1-\frac{4p\gamma(\delta_1-\delta_0)}{(\delta_0-\delta_1-\gamma)^2}}\right)\\
&\quad \ge \frac{1}{2}\left(\delta_0+\delta_1+\gamma-|\delta_0-\delta_1-\gamma|+\frac{2p\gamma(\delta_1-\delta_0)}{|\delta_0-\delta_1-\gamma|}\right) = \delta_0+\frac{p\gamma(\delta_1-\delta_0)}{|\delta_0-\delta_1-\gamma|}.
\end{aligned}$$

Therefore, for every $p>0$, we can choose $\delta_0 < \delta_{c2}$ sufficiently close to $\delta_{c2}$ so that $\lambda_{\max,\mu}(\delta_0,\delta_1,\gamma,p) > \delta_{c2}$. Therefore, as above, the process $\Psi^{\gamma,p,A}_{\delta_0,\delta_1}$ dies out and therefore

$$\lim_{\delta_0\uparrow\delta_{c2}} p_{c2}(\delta_0,\delta_1,\gamma) < p.$$

Since $p>0$ was arbitrary, $\lim_{\delta_0\uparrow\delta_{c2}} p_{c2}(\delta_0,\delta_1,\gamma)=0$ and we are done. $\square$

**6. Open questions.** We here list some open questions related to the results of this paper.

1. Do either of the critical values $p_{c1}$ and $p_{c2}$ depend on the initial state of the background process?
2. Instead of studying the CPREE model one could study other interacting particle systems such as a stochastic Ising model in a random evolving environment.
3. Is it possible to generalize the model used for the background process in some way? For instance, can we analyze the situation where we allow more than 2 different states?

**Acknowledgments.** I thank Jeff Steif for all the suggestions he had during the work on this paper as well as for patiently reading this manuscript many times before its completion. I would also like to thank the referee for a careful reading as well as for giving many suggestions and comments, for instance, Examples 5.1 and 5.2 as well as shortening the proof of Theorem 1.8.

Department of Mathematics
Chalmers University of Technology
S-412 96 Göteborg
Sweden
E-mail: broman@math.chalmers.se
URL: http://www.math.chalmers.se/~broman/